\documentclass[12pt, reqno]{amsart}
\usepackage{amssymb, amsthm, amsmath, amsfonts}
\usepackage{hyperref,verbatim}
\usepackage{array, epsfig}
\usepackage{bbm}

\usepackage{color}

  \evensidemargin
\oddsidemargin
\begin{document}

\newcommand{\be}{\begin{equation}}
\newcommand{\ee}{\end{equation}}
\newcommand{\bea}{\begin{eqnarray}}
\newcommand{\eea}{\end{eqnarray}}
\newcommand{\beaa}{\begin{eqnarray*}}
\newcommand{\eeaa}{\end{eqnarray*}}

\renewcommand{\proofname}{\bf Proof}
\newtheorem*{rem*}{Remark}
\newtheorem*{cor*}{Corollary}
\newtheorem{cor}{Corollary}
\newtheorem*{conj}{Conjecture}
\newtheorem{prop}{Proposition}
\newtheorem*{prop1}{Proposition 1}
\newtheorem*{def*}{Definition}
\newtheorem{lem}{Lemma}
\newtheorem*{lem1'}{Lemma $\mathbf{1^\prime}$}
\newtheorem{theo}{Theorem}
\newfont{\zapf}{pzcmi}

\newcommand{\Var}{\mathop{\mathrm{Var}}\nolimits}
\newcommand{\supp}{\mathop{\mathrm{supp}}\nolimits}
\newcommand{\Cov}{\mathop{\mathrm{Cov}}\nolimits}
\newcommand{\sgn}{\mathop{\mathrm{sgn}}\nolimits}
\newcommand{\law}{\mathop{\mathrm{Law}}\nolimits}

\def\R{\mathbb{R}}
\def\Z{\mathbb{Z}}
\def\N{\mathbb{N}}
\def\E{\mathbb{E}}
\def\P{\mathbb{P}}
\def\V{\mathbb{D}}
\def\S{\mathbb{S}}
\def\I{\mathbbm{1}}
\newcommand{\D}{\hbox{\zapf D}}
\newcommand{\Pt}{\widetilde{\mathbb{P}}}
\newcommand{\Et}{\widetilde{\mathbb{E}}}

\newcommand{\bt}{\begin{theo}}
\newcommand{\et}{\end{theo}}
\newcommand{\bl}{\begin{lem}}
\newcommand{\el}{\end{lem}}
\newcommand{\bc}{\begin{cor*}}
\newcommand{\ec}{\end{cor*}}
\newcommand{\br}{\begin{rem*}}
\newcommand{\er}{\end{rem*}}
\newcommand{\bp}{\begin{proof}}
\newcommand{\ep}{\end{proof}}
\newcommand{\bes}{\begin{ex}}
\newcommand{\ees}{\end{ex}}

\title{On the weak limit law of the maximal uniform $k$-spacing}
\author[A. Mijatovi\'c and V. Vysotsky]{Aleksandar Mijatovi\'c$^\dagger$ and Vladislav Vysotsky$^\ddagger$}
\address{Arizona State University$^\ddagger$, Imperial College London$^{\dagger, \ddagger}$, St.Petersburg Department of Steklov Mathematical Institute$^{\ddagger}$}
\thanks{AM is supported in part by the Humboldt Foundation Research Fellowship GRO/1151787 STP.
The work of VV is supported by Marie Curie IIF Grant 628803 by European Commission and supported in part by Grant 13-01-00256 by RFBR}

\begin{abstract}
This paper gives a simple proof of a limit theorem for the length of the largest interval straddling a fixed number of i.i.d. points uniformly distributed on a unit interval. The key step in our argument is a classical theorem of~\cite{Watson} on the maxima of $m$-dependent stationary stochastic sequences.
\end{abstract}

\maketitle

\section{Introduction and the main result}
Both the distributional and asymptotic theories of spacings between consecutive order statistics of a sample
of i.i.d. random variables play a central role in classical probability theory and mathematical statistics, see~\cite{Pyke},~\cite[Sec.~18--21]{SW} and the references therein. A deep understanding of this subject has been achieved over the past decades. In particular,~\cite{Devroye, Deheuvels} give a very fine description of the almost sure behaviour (as the sample size tends to infinity) 
of the maximal spacing between the ordered statistics of uniform random variables. The laws of iterated logarithms proved in these papers for the maximal spacings are further extended in~\cite{DD} to analogous statements on the maximum of $k$ consecutive spacings (called $k$-spacings).


In this note we prove a weak limit theorem for the maximal $k$-spacings. To the best of our knowledge, no result of this type was available in the past; it is truly surprising that this problem was not even mentioned in~\cite{DD}.

More precisely, let $U_1, \dots, U_n$ be i.i.d. random variables that are uniformly distributed on $[0,1]$. Denote by
$$U_{1:n} \le \dots \le U_{n:n}$$ 
their order statistics, which are the elements of $U_1, \dots, U_n$ arranged in the ascending order, and define 
$U_{0:n}:=0$,
$U_{(n+1):n}:=1$.
The maximal spacing
$M^{(1)}_n:= \max_{0 \le i \le n} (U_{(i+1):n} - U_{i:n})$
is the lenght of the longest interval containing no points of the sample 
$U_1, \dots, U_n$.
The classical representation of the uniform spacings given below in~\eqref{eq:representation}, 
which relates $M^{(1)}_n$ to the maximum of i.i.d. exponential random
variables, together with the law of large numbers easily yields 
$$n M^{(1)}_n - \log n
\stackrel{d}{\longrightarrow} G,$$
where 
$G$
follows a standard Gumbel  distribution $\P(G \le x)= \exp\left(-e^{-x}\right), x \in \R$. 

We study an analogous weak limit of the \textit{maximal $k$-spacing}, that is the length of the largest open subinterval of $[0,1]$ that contains $k-1$ uniform points:
$$M^{(k)}_n:= \max_{0 \le i \le n+1-k} (U_{(i+k):n} - U_{(i):n}).$$
Our main result is as follows. 
\begin{theo} \label{Thm main}
Let $G$ be a random variable that follows a standard Gumbel distribution. For any 
integer $k \ge 1$, it holds 
$$n M^{(k)}_n - \log n - (k-1)\log\log n +\log (k-1)! \stackrel{d}{\longrightarrow} G\qquad\text{as $n\to\infty$.}$$
\end{theo}

We will use the following well-known fact: the uniform spacings are represented as 
\be
\Bigl(U_{1:n} - U_{0:n}, \dots, U_{n:n} - U_{(n-1):n}\Bigr) \stackrel{d}{=} \Bigl ( \frac{X_1}{X_1+\dots+X_{n+1}}, \dots, \frac{X_n}{X_1+\dots+X_{n+1}}\Bigr),
\label{eq:representation}
\ee 
where $X_1, X_2, \dots$ are i.i.d. standard exponential random variables; 
moreover, the random vector on the right-hand side is independent of the sum 
$X_1+\dots+X_{n+1}$, see e.g.~\cite[Sec.~4.1]{Pyke}.


To discuss the statement of Theorem~\ref{Thm main}, consider the simplest case that $k=2$ on the largest interval straddling a single uniform point. It is not hard to show that
$A_n:= \max_{1 \le i \le n} (X_{2i-1} + X_{2i})$, $B_n:= \max_{1 \le i \le n} (X_{2i} + X_{2i+1})$, which 
are maxima of i.i.d. gamma random variables, satisfy
$A_n - \log n - \log \log n \stackrel{d}{\longrightarrow} G$ and $B_n - \log n - \log \log n \stackrel{d}{\longrightarrow} G$.

The crucial observation is that $A_n - \log n - \log \log n$ and $B_n - \log n - \log \log n$ are {\it asymptotically independent}.
Then
\be 
\nonumber
M^{(2)}_n=\frac{\max(A_{\lfloor n/2 \rfloor}, B_{\lfloor (n-1)/2 \rfloor})}{X_1+\dots+X_{n+1}},
\ee 
and hence the law of large numbers and the continuous mapping theorem imply\footnote{See the analogous argument after \eqref{eq: max Y} below.} 
$$nM^{(2)}_n - \log (n/2) - \log \log(n/2) \stackrel{d}{\longrightarrow} \max(G_1, G_2),$$ 
where $G_1$ and $G_2$ are i.i.d. random variables with a standard Gumbel distribution. Since $\max(G_1, G_2) \stackrel{d}{=} \log 2 + G$, 
Theorem~\ref{Thm main} follows in the case that $k=2$.

The asymptotic independence
of
$A_n$
and
$B_n$
is non-trivial and somewhat unexpected. Our initial approach to the proof of Theorem~\ref{Thm main} rested on establishing this
property using the specific structure of these random variables. However, once the classical result~\cite{Watson}
on the maxima of $m$-dependent stationary sequences came to our attention, we understood that our
Theorem~\ref{Thm main} can be established as a direct consequence\footnote{The asymptotic independence appears not to be an easy consequence of the result in~\cite{Watson}.}. We describe this shorter and easier proof in the next section.

\section{Proofs}
\label{sec:proofs}

We start by recalling the result from~\cite{Watson}.
Random variables 
$Y_1, Y_2, \dots$
are said to be $m$-\textit{dependent}
if 
$|i-j|>m$
implies that $Y_i$
and
$Y_j$
are independent.

\begin{theo}
\label{thm:Watson}
For any 
$m\ge 1$,
let $Y_1, Y_2, \dots$ 
be 
a strictly stationary sequence of $m$-dependent unbounded random variables. Assume that
\begin{equation}
\label{eq:As1}
\lim_{y\to\infty} \max_{1\leq|i-j|\leq m}\P(Y_j>y|Y_i>y)=0.
\end{equation}
Then for any positive numbers $\xi, y_1, y_2, \dots$ satisfying
\begin{equation}
\label{eq:As2}
\lim_{n \to \infty} n\P(Y_1>y_n) = \xi,
\end{equation}
it holds 
$$
\lim_{n\to\infty}\P\left(\max_{1\leq i \leq n}Y_i\leq y_n\right) = \exp(-\xi).
$$
\end{theo}

The theorem says that the maximum of $m$-dependent stationary random variables
has the same weak limit as the maximum of an i.i.d. sequence with
the same common distribution. Although the actual theorem of~\cite{Watson}
makes a more restrictive assumption $\xi = n\P(Y_1>y_n)$ for all $n \ge 1$,
which may even be impossible to satisfy for certain $\xi$, the presented version easily
follows by the monotonicity of distribution functions and the continuity of
$\exp(-\xi)$. 

\medskip
The aim is to apply Theorem~\ref{thm:Watson} to the $(k-1)$-dependent stationary sequence of moving sums
\begin{equation}
\label{eq:Def_Y}
Y_i:=\sum_{\ell=i}^{i+k-1}X_\ell, \qquad i \ge 1,
\end{equation}
and the numbers 
\begin{equation}
\label{eq:seq_cn}
\xi:=e^{-x}, \qquad y_n:=\log n + (k-1)\log\log n - \log (k-1)! + x
\end{equation}
for any fixed real $x$.

Note first that $Y_i$ are gamma random variables with densities $f_k$, where $f_\theta(y):=y^{\theta-1}e^{-y}/\Gamma(\theta)$ 
for any positive 
$y$ and $\theta$. Then it is straightforward to check using L'Hopital's rule that
\be
\label{eq: tail}
\P(Y_1 > y) \sim \frac{y^{k-1} e^{-y}}{(k-1)!}, \qquad y \to \infty
\ee
(where by $\sim$ we mean that the ratio tends to $1$),
hence \eqref{eq:As2} holds by
$$
\lim_{n \to \infty} n\P\left(Y_1>y_n\right) = \lim_{n \to \infty} \frac{n y_n^{k-1}e^{-y_n}}{(k-1)!} = e^{-x} \lim_{n \to \infty} \left(\frac{y_n}{\log n}\right)^{k-1} = \xi.
$$

It remains to check the assumption \eqref{eq:As1}. For any integer $1 \le a \le k-1$, we have
$$Y_{a+1} = Y_1\Bigl( 1 - \frac{X_1 + \dots + X_a}{X_1+ \dots + X_k} \Bigr) + (X_{k+1} + \dots + X_{k+a}).$$
Hence
$$(Y_1, Y_{a+1}) \stackrel{d}{=} \bigl( Y_1, Y_1 (1 - U_{a:(k-1)}) + Z_a)\bigr),$$
where the three random variables in the r.h.s. are mutually independent and $Z_a$ has a gamma distribution with density $f_a$. By \eqref{eq: tail}, for any $\varepsilon > 0$ there exists an $R>0$ such that 
$$\P(Y_1>y+R) \leq \varepsilon \P(Y_1>y )   \qquad \mbox{for all $y$ large enough}.$$ 
Then \eqref{eq:As1} follows as for such $y$,
\beaa
\P(Y_i > y, Y_{a+i} >y) &\le& \P(y<Y_1\leq y+R, Y_{a+1}>y) +  \P(Y_1 > y+R)\\
& \le & 
 \int_y^{y+R} \P \bigl(Z_a> y-x(1- U_{a:(k-1)})\bigr) f_k(x) dx + \varepsilon \P(Y_1>y) \\
&\le& \bigl( \P\bigl(Z_{k-1}> y U_{1:(k-1)} - R\bigr)+\varepsilon\bigr)\cdot \P(Y_1>y).
\eeaa

Thus we showed that Theorem~\ref{thm:Watson} applies to the sequence 
$Y_1, Y_2, \dots$ 
defined 
in~\eqref{eq:Def_Y}, hence combined with \eqref{eq:seq_cn} this implies
\be
\label{eq: max Y}
\max \limits_{1 \le i \le n+1-k} Y_i - \log n - (k-1)\log\log n + \log (k-1)! \stackrel{d}{\longrightarrow} G.
\ee
Then by~\eqref{eq:representation}, we find
$$nM^{(k)}_n= \frac{n}{X_1+\cdots+X_{n+1}} \max \limits_{1 \le i \le n+1-k} Y_i.$$
Now Theorem~\ref{Thm main} follows by \eqref{eq: max Y}, the law of large numbers, the continuous mapping theorem, and the relation
$$\log n \Bigl ( \frac{n}{X_1 + \dots +X_{n+1} } - 1 \Bigr) = \frac{\log n}{\sqrt{n}} \cdot \frac{(n- ( X_1 + \dots +X_{n+1}))/\sqrt{n} }{(X_1 + \dots +X_{n+1})/n} \stackrel{d}{\longrightarrow} 0,$$ which itself holds by
the law of large numbers, and the central limit theorem.

\bibliographystyle{apalike}
\bibliography{biblio}

\end{document}